\documentstyle[11pt,amssymb]{article}
\newcommand{\beql}[1]{\begin{equation}\label{#1}}
\newcommand{\eeq}{\end{equation}}
\newcommand{\comment}[1]{}
\newcommand{\Ds}{\displaystyle}
\newcommand{\eqref}[1]{{\rm (\ref{#1})}}

\newcommand{\Abs}[1]{{\left|{#1}\right|}}
\newcommand{\Lone}[1]{{\left\|{#1}\right\|_{L^1}}}

\newcommand{\Qed}{\ \\\mbox{$\blacksquare$}}

\newcommand{\Set}[1]{{\left\{{#1}\right\}}}

\newcommand{\ToAppear}[1]{\raisebox{15mm}[10pt][0mm]{\makebox[0mm]%
	{\makebox[\textwidth][r]{\small #1}}}}

\newcommand{\RR}{{\Bbb R}}
\newcommand{\CC}{{\Bbb C}}
\newcommand{\ZZ}{{\Bbb Z}}

\newcommand{\TT}{{\Bbb T}}

\newcommand{\one}{{\bf 1}}

\newcommand{\inner}[2]{{\langle #1, #2 \rangle}}
\newcommand{\Inner}[2]{{\left\langle #1, #2 \right\rangle}}
\newcommand{\dens}{{\rm dens\,}}

\newcommand{\supp}{{\rm supp\,}}

\newcommand{\ft}[1]{\widehat{#1}}
\newcommand{\FT}[1]{\left(#1\right)^\wedge}

\newcounter{open}
\newcounter{dfn}
\def\thedfn{\arabic{dfn}}
\newenvironment{dfn}{
  \sf
  \vskip 0.10in
  \refstepcounter{dfn}
  \noindent{\bf Definition \thedfn \ }
}{\vskip 0.10in}

\newcounter{obs}
\def\theobs{\arabic{obs}}

\newcounter{thm}
\newenvironment{thm}{
  \sf
  \vskip 0.10in
  \refstepcounter{thm}
  \noindent{\bf Theorem\ }
}{\vskip 0.10in}

\newcounter{mysec}
\def\themysec{\arabic{mysec}}
\newcommand{\mysection}[1]{
  \vskip 0.25in
  \refstepcounter{mysec}\centerline{\large\bf \S\themysec.\ {#1}}\par
  \nopagebreak
  \addcontentsline{toc}{section}{{\bf \themysec.}\ {#1}}
}

\newcounter{mysubsec}[mysec]
\newtheorem{theorem}{Theorem}
\newtheorem{corollary}{Corollary}

\newtheorem{conj}{Conjecture}
\begin{document}

\begin{center}
{\Large \bf Packing, tiling, ortho\ToAppear{
{\tt Bull.\ LMS, to appear}}gonality and completeness}\\
\ \\
{\sc Mihail N. Kolountzakis\footnote{Partially supported%
by the U.S. National Science Foundation,
under grant DMS 97-05775.
Most of this work was carried out while the author
was visiting the Univ.\ of Illinois at Urbana-Champaign,
in Fall 1998-99.}}\\
Department of Mathematics,\\
University of Crete,\\
Knossos Ave., 714 09 Iraklio,\\
Greece.\\
\ \\
E-mail: {\tt kolount@math.uch.gr}\\
\ \\
\small March 1999; revised August 1999
\end{center}

\begin{abstract}
Let $\Omega \subseteq \RR^d$ be an open set of measure $1$. An open
set $D \subseteq \RR^d$ is called a ``tight orthogonal packing region''
for $\Omega$ if $D-D$ does not intersect the zeros of the Fourier
Transform of the indicator function of $\Omega$ and
$D$ has measure $1$.
Suppose that $\Lambda$ is a discrete subset of $\RR^d$.
The main contribution of this paper is a new way of proving
the following result (proved by different methods
by Lagarias, Reeds and Wang \cite{LRW98} and,
in the case of $\Omega$ being the cube, by Iosevich
and Pedersen \cite{IP98}): $D$ tiles $\RR^d$ when translated at the locations
$\Lambda$ if and only if the set of exponentials
$E_\Lambda = \Set{\exp 2\pi i \inner{\lambda}{x}:\ \lambda\in\Lambda}$
is an orthonormal basis for $L^2(\Omega)$.
(When $\Omega$ is the unit cube in $\RR^d$ then it is
a tight orthogonal packing region of itself.)
In our approach orthogonality of $E_\Lambda$ is viewed
as a statement about ``packing'' $\RR^d$ with translates of a
certain nonnegative function and,
additionally, we have completeness of $E_\Lambda$ in $L^2(\Omega)$
if and only if the above-mentioned packing is
in fact a tiling.
We then formulate the tiling condition in Fourier Analytic language
and use this to prove our result.
\end{abstract}

\mysection{Introduction}\label{sec:intro}

\noindent
{\bf Notation.}\\
Let $\Omega \subset \RR^d$ be measurable of measure $1$.
The Hilbert space $L^2(\Omega)$ is equipped with the inner
product
$$
\inner{f}{g}_\Omega = \int_\Omega f(x) \overline{g(x)}~dx.
$$
Define
$$
e_\lambda(x) = \exp 2\pi i \inner{\lambda}{x},
$$
and, for $\Lambda \subseteq \RR^d$,
$$
E_\Lambda = \Set{e_\lambda:\ \lambda\in\Lambda}.
$$
For every continuous function $h:\RR^d\to\CC$ we write
$$
Z(h) = \Set{x\in\RR^d:\ h(x) = 0}.
$$
Whenever we fail to mention it it should be understood that
the measure of $\Omega \subset \RR^d$ is equal to $1$.

The indicator function of a set $E$ is denoted by $\one_E$.

We denote by $B_r(x)$ the ball in $\RR^d$ of radius $r$ centered at $x$.

When $A$ and $B$ are two sets in $\RR^d$ we write $A+B$ for the
set of all sums $a+b$, $a\in A$, $b\in B$. Similarly we write
$A-B$ for all differences $a-b$, $a\in A$, $b\in B$.
For $\lambda \in \RR$ we denote by $\lambda A$ the set
$\Set{\lambda a:\ a\in A}$.

If $O$ is an open set in $\RR^d$ we denote by $C_c^\infty(O)$
the set of all infinitely differentiable functions with support
contained in $O$.

%%%%%%%%%%%%%%%%%%%%%%%%%%%%%%%%%%%%%%%%%%%%%%%%%%%%%%%%%%%%%%%%
%%%%%%%%%%%%%%% dfn SPECTRAL SET
%%%%%%%%%%%%%%%%%%%%%%%%%%%%%%%%%%%%%%%%%%%%%%%%%%%%%%%%%%%%%%%%
\begin{dfn}\label{def:spectral-set}
(Spectral sets)\\
Suppose that $\Omega$ is a measurable set of measure $1$.
We call $\Omega$ \underline{spectral} if $L^2(\Omega)$ 
has an orthonormal basis
$E_\Lambda=\Set{e_\lambda:\ \lambda\in\Lambda}$ of exponentials.
The set $\Lambda$ is then called a \underline{spectrum} for $\Omega$.
\end{dfn}
We can always restrict our attention to sets $\Lambda$ containing $0$
and we shall do so without further mention.

%%%%%%%%%%%%%%%%%%%%%%%%%%%%%%%%%%%%%%%%%%%%%%%%%%%%%%%%%%%%%%%%
%%%%%%%%%%%%%%% dfn PACKING, TILING
%%%%%%%%%%%%%%%%%%%%%%%%%%%%%%%%%%%%%%%%%%%%%%%%%%%%%%%%%%%%%%%%
\begin{dfn}\label{def:packing-tiling}
(Packing and tiling by nonnegative functions)\\
(i) A nonnegative measurable function $f: \RR^d \to \RR^+$
(the ``tile'') is said to \underline{pack} a region $S \subseteq \RR^d$
with the set $\Lambda \subseteq \RR^d$ (the ``set of translates'') if
$$
\sum_{\lambda \in \Lambda} f(x-\lambda) \le 1 \ \mbox{for a.e.\ $x\in S$}.
$$
In this case we write ``$f + \Lambda$ packs $S$''. When $S$ is omitted
we understand $S=\RR^d$.\\
(ii) A nonnegative measurable function $f: \RR^d \to \RR^+$ is said
to \underline{tile} a region $S \subseteq \RR^d$ at \underline{level}
$\ell$ with the set $\Lambda \subseteq \RR^d$ if
$$
\sum_{\lambda \in \Lambda} f(x-\lambda) = \ell \ \mbox{for a.e. $x\in S$}.
$$
(When not specified $\ell = 1$.)
Again we write ``$f + \Lambda$ tiles $S$ at level $\ell$'' (or
$f + \Lambda = \ell S$) and $S = \RR^d$ is
understood when $S$ is omitted.

If $f = \one_E$ is the indicator function of a measurable set $E$
we also say that ``$E+\Lambda$'' is a packing (resp. tiling) instead
of ``$\one_E + \Lambda$ is a packing'' (resp. tiling).
\end{dfn}

%%%%%%%%%%%%%%%%%%%%%%%%%%%%%%%%%%%%%%%%%%%%%%%%%%%%%%%%%%%%%%%%
%%%%%%%%%%%%%%% conj FUGLEDE
%%%%%%%%%%%%%%%%%%%%%%%%%%%%%%%%%%%%%%%%%%%%%%%%%%%%%%%%%%%%%%%%
The following conjecture of Fuglede \cite{F74} is still
unresolved and has provided the motivation of the problem we
deal with in this paper.
\begin{conj}\label{conj:fuglede}
{\sf (Fuglede)}\\
Let $\Omega$ be a bounded open set of measure $1$.
Then $\Omega$ is spectral if and only if $\Omega$
tiles $\RR^d$ by translation.
\end{conj}
As an example of a spectral set in $\RR^2$ we give the open
unit square $(-1/2, 1/2)^2$ which tiles the plane when translated by
$\ZZ^2$ and also has $\ZZ^2$ as its spectrum.
(Note that in Fuglede's conjecture the set of translates
by which a tile $\Omega$ tiles space need not be the same as its spectrum.)
Fuglede \cite{F74} proved that the triangle and the disk in the plane
are both not spectral.
Further in this direction of confirming the conjecture we mention
that Iosevich, Katz and Pedersen \cite{IKP99} have recently proved
that the ball in $\RR^d$ is not spectral and the author \cite{nonsym}
proved that any non-symmetric convex domain in $\RR^d$ is not spectral.
Convex domains that tile space by translation are known to be
necessarily symmetric (see \cite{M80}).

A related problem is, given a specific
set $\Omega$ that tiles space by translation, to determine its spectra.
Because of its simplicity the cube has been studied the most.
Lagarias, Reeds and Wang \cite{LRW98} and Iosevich and Pedersen \cite{IP98}
recently proved that if $Q = (-1/2,1/2)^d$ is the unit cube in $\RR^d$
then $Q+\Lambda$ is a tiling if and only if $E_\Lambda$ is an orthonormal
basis for $Q$.
(We remark here that there exist ``exotic'' translational tilings by
the unit cube which are non-lattice, see \cite{LS}.)
This had been conjectured by Jorgensen and Pedersen \cite{JP}
where it was proved for dimension $d\le 3$.
The purpose of our paper is to give an alternative and, perhaps, more
illuminating proof of this fact, which is based on a characterization
of translational tiling by a Fourier Analytic criterion.

We follow the terminology of \cite{LRW98}.

The two basic tools in this paper are Theorem \ref{th:tiling-from-packing}
and Theorem \ref{th:supp-equiv-tiling}.
Theorem \ref{th:tiling-from-packing}, which is interesting and rather
unexpected in itself, concerns tilings by two different tiles and states
that if two tiles $A$ and $B$, of the same volume, both pack with $\Lambda$
then they either both tile or both do not.
This is almost obvious if $\Lambda$ is periodic but, interestingly, it holds in
general and its proof is rather simple.
Theorem \ref{th:supp-equiv-tiling} is the final form for our
Fourier Analytic criterion.

%%%%%%%%%%%%%%%%%%%%%%%%%%%%%%%%%%%%%%%%%%%%%%%%%%%%%%%%%%%%%%%%
%%%%%%%%%%%%%%% dfn DENSITY
%%%%%%%%%%%%%%%%%%%%%%%%%%%%%%%%%%%%%%%%%%%%%%%%%%%%%%%%%%%%%%%%
\begin{dfn}
(Density)\\
A set $\Lambda\subseteq\RR^d$ has \underline{asymptotic} density
$\rho$ if
$$
\lim_{R\to\infty} {\#(\Lambda \cap B_R(x)) \over \Abs{B_R(x)}} \to
\rho
$$
uniformly in $x\in\RR^d$.\\
We say that $\Lambda$ has \underline{(uniformly) bounded density} if the
fraction
above is bounded by a constant $\rho$ uniformly for $x\in\RR$ and
$R>1$.
We say then that $\Lambda$ has density (uniformly) bounded by $\rho$.
\end{dfn}
If $f \ge 0$ then it is clear that if $f + \Lambda$ is a tiling at level
$\ell>0$
then
$\Lambda$ has asymptotic density equal to $\ell / \int f$.

%%%%%%%%%%%%%%%%%%%%%%%%%%%%%%%%%%%%%%%%%%%%%%%%%%%%%%%%%%%%%%%%
%%%%%%%%%%%%%%% section CONNECTIONS
%%%%%%%%%%%%%%%%%%%%%%%%%%%%%%%%%%%%%%%%%%%%%%%%%%%%%%%%%%%%%%%%
\mysection{Packing, tiling, orthogonality and completeness}
\label{sec:connections}

Let $\Omega$ be a measurable set in $\RR^d$ of measure $1$
and $\Lambda$ be a countable set of points in $\RR^d$.

The set $E_\Lambda = \Set{e_\lambda(x):\ \ \lambda\in\Lambda}$
is an orthogonal set of exponentials for $\Omega$ if and only if
$$
\sum_{\lambda\in\Lambda} \Abs{\inner{e_x}{e_\lambda}_\Omega}^2 \le 1,
$$
for each $x \in \RR^d$.
Since
$$
\inner{e_x}{e_\lambda}_\Omega = 
  \int e^{2\pi i(x-\lambda)t} \one_\Omega(t)~dt =
  \ft{\one_\Omega}(\lambda - x),
$$
we conclude that
$\Lambda$ is an orthogonal set for $\Omega$
if and only if $\Abs{\ft{\one_\Omega}}^2 + \Lambda$ is a packing
of $\RR^d$.

In this case $\Lambda$ has uniformly bounded density.

Similarly, $\Lambda$ is a spectrum of $\Omega$
($E_\Lambda$ is orthogonal and complete) if and only if
$$
\sum_{\lambda\in\Lambda} \Abs{\inner{e_x}{e_\lambda}_\Omega}^2 = 1,
$$
for all $x \in \RR^d$.
That is, $\Lambda$ is a spectrum of $\Omega$ if and only if
$\Abs{\ft{\one_\Omega}}^2 + \Lambda$ is a tiling of $\RR^d$.

\begin{dfn}\label{def:packing-region}
The open set $D$ is called an \underline{orthogonal packing region} for
$\Omega$
if
$$
(D-D) \cap Z(\ft{\one_\Omega}) = \emptyset.
$$
\end{dfn}

By the definition of an orthogonal packing region $D$ for $\Omega$,
if $\Lambda$ is an orthogonal set of exponentials for $\Omega$ then
$D + \Lambda$ is a packing of $\RR^d$.
Indeed, if $\lambda, \mu \in \Lambda$, $\lambda \neq \mu$, then
$\ft{\one_\Omega}(\lambda-\mu) = 0$,
since $\Abs{\ft{\one_\Omega}}^2 + \Lambda$ is a packing,
which implies $\lambda - \mu \in Z(\ft{\one_\Omega})$ which is
disjoint from $D-D$.
Hence $(\Lambda-\Lambda)\cap(D-D) = \Set{0}$
and $D+\Lambda$ is a packing.

We summarize these observations in the following theorem.
\begin{theorem}\label{th:connections}
Let $\Omega$ be a measurable set in $\RR^d$ of
measure $1$ and $\Lambda \subset \RR^d$
be countable.
\begin{enumerate}
\item
$E_\Lambda$ is an orthogonal set for $L^2(\Omega)$ if and only if
$\Abs{\ft{\one_\Omega}}^2 + \Lambda$ is a packing.
\item
$E_\Lambda$ is an orthonormal basis for $L^2(\Omega)$ (a spectrum
for $\Omega$) if and only if
$\Abs{\ft{\one_\Omega}}^2 + \Lambda$ is a tiling.
\item
If $D$ is an orthogonal packing region for $\Omega$ and $E_\Lambda$
is an orthogonal set in $L^2(\Omega)$ then
$D + \Lambda$ is a packing.
\end{enumerate}
\end{theorem}

%%%%%%%%%%%%%%%%%%%%%%%%%%%%%%%%%%%%%%%%%%%%%%%%%%%%%%%%%%%%%%%%
%%%%%%%%%%%%%%% section RESULT ABOUT PACKING AND TILING
%%%%%%%%%%%%%%%%%%%%%%%%%%%%%%%%%%%%%%%%%%%%%%%%%%%%%%%%%%%%%%%%

\mysection{A result about packing and tiling by two
different tiles}
\label{sec:interesting}

The following theorem is a crucial tool for the results of this
paper.
It is intuitively clear when $\Lambda$ is a periodic set but it is,
perhaps, suprising that it holds without any assumptions on
the set $\Lambda$.
Its proof is very simple.

%%%%%%%%%%%%%%%%%%%%%%%%%%%%%%%%%%%%%%%%%%%%%%%%%%%%%%%%%%%%%%%%
%%%%%%%%%%%%%%% thm about PACKING & TILING, two functions
%%%%%%%%%%%%%%%%%%%%%%%%%%%%%%%%%%%%%%%%%%%%%%%%%%%%%%%%%%%%%%%%
\begin{theorem}\label{th:tiling-from-packing}
If $f, g \ge 0$, $\int f(x) dx = \int g(x) dx = 1$ and 
both $f+\Lambda$ and $g+\Lambda$ are packings of $\RR^d$, then
$f+\Lambda$ is a tiling if and only if $g+\Lambda$ is a tiling.
\end{theorem}
{\bf Proof.}
We first show that, under the assumptions of the Theorem,
\beql{supp-tiling}
\mbox{$f+\Lambda$ tiles $-\supp g$ } \Longrightarrow
\mbox{ $g+\Lambda$ tiles $-\supp f$}.
\eeq
Indeed, if $f+\Lambda$ tiles $-\supp g$ then
$$
1 = \int g(-x) \sum_{\lambda\in\Lambda} f(x-\lambda) ~dx = 
 \sum_{\lambda\in\Lambda} \int g(-x) f(x-\lambda) ~dx,
$$
which, after the change of variable $y = -x+\lambda$, gives
$$
1 = \int f(-y) \sum_{\lambda\in\Lambda} g(y-\lambda) ~dy.
$$
This in turn implies,
since $\sum_{\lambda\in\Lambda}g(y-\lambda) \le 1$,
that $\sum_\lambda g(y-\lambda) = 1$ for
a.e.\ $y \in -\supp f$.

To complete the proof of the theorem, notice that if $f+\Lambda$ is a tiling
of $\RR^d$ and $a \in \RR^d$ is arbitrary then
both $f(x-a) + \Lambda$ and $g(x-a) + \Lambda$ are packings and
$f+\Lambda$ tiles $-\supp g(x-a) = -\supp g - a$.
We conclude that $g(x-a) + \Lambda$ tiles $-\supp f$, or
$g + \Lambda$ tiles $-\supp f - a$. Since $a\in\RR^d$ is arbitrary
we conclude that $g + \Lambda$ tiles $\RR^d$.
\Qed

%%%%%%%%%%%%%%%%%%%%%%%%%%%%%%%%%%%%%%%%%%%%%%%%%%%%%%%%%%%%%%%%
%%%%%%%%%%%%%%% section FOURIER ANALYTIC CRITERIA
%%%%%%%%%%%%%%%%%%%%%%%%%%%%%%%%%%%%%%%%%%%%%%%%%%%%%%%%%%%%%%%%

\mysection{Fourier Analytic criteria for tiling}
\label{sec:fourier}

The action of a tempered distribution (see \cite{R73})
$\alpha$ on a Schwartz function $\phi$
is denoted by $\alpha(\phi)$. The Fourier Transform of $\alpha$ is
defined by the equation
$$
\ft\alpha(\phi) = \alpha(\ft\phi).
$$
The support $\supp\alpha$ is the smallest closed set $F$ such that for any
smooth $\phi$ of compact support contained in the open set $F^c$ we have
$\alpha(\phi)=0$.

If $f \in L^1(\RR^d)$ and $\widehat{f}$ is $C^\infty$ then,
if $\Lambda\subset\RR^d$ is a discrete set, the following theorem,
first proved in \cite{KL96} in dimension $1$, gives necessary
and sufficient conditions for $f+\Lambda$ to be a tiling.
We give the proof here for completeness.
%%%%%%%%%%%%%%%%%%%%%%%%%%%%%%%%%%%%%%%%%%%%%%%%%%%%%%%%%%%%%%%%
%%%%%%%%%%%%%%% thm KL96
%%%%%%%%%%%%%%%%%%%%%%%%%%%%%%%%%%%%%%%%%%%%%%%%%%%%%%%%%%%%%%%%
\begin{theorem}\label{th:tiling-smooth}
Suppose $f\in L^1(\RR^d)$ and $\widehat{f} \in C^\infty$.
Suppose also that $\Lambda$ is a discrete subset of $\RR^d$
of bounded density.
Write $\delta_\Lambda$ for the tempered distribution
$\sum_{\lambda\in\Lambda} \delta_\lambda$ and $\widehat{\delta_\Lambda}$
for its Fourier Transform.\\
(i) If $f+\Lambda$ is a tiling then
\beql{sp-cond-1}
\supp \widehat{\delta_\Lambda} \subseteq \Set{\widehat f = 0} \cup \Set{0}.
\eeq
(ii) If $\widehat{\delta_\Lambda}$ is locally a measure then
\eqref{sp-cond-1} implies that $f+\Lambda$ is a tiling.
\end{theorem}
Notice that whenever $f$ has compact support
the function $\widehat{f}$ is smooth.

\noindent
{\bf Proof of Theorem \ref{th:tiling-smooth}.}
(i) If $f + \Lambda$ is a tiling then $1 = f*\delta_\Lambda$, hence,
taking Fourier Transforms, $\delta_0 = \widehat{f}\widehat{\delta_\Lambda}$.
Take $\phi \in C_c^\infty(\RR^d\setminus K)$, where
$$
K = \Set{\widehat f = 0} \cup \Set{0}.
$$
Then
$$
\widehat{\delta_\Lambda}(\phi) = 
(\widehat{f}\widehat{\delta_\Lambda}) \left({\phi\over\widehat{f}}\right) =
\delta_0\left({\phi\over\widehat{f}}\right) = {\phi\over\widehat{f}}(0) = 0.
$$
This proves \eqref{sp-cond-1}. Note that it was crucial in the proof
that $\Ds{\phi\over\widehat{f}}$ is a function in 
$C_c^\infty(\RR^d\setminus K)$ and this is why we demanded that
$\widehat{f}$ is smooth.

(ii) Take $\phi$ to be a Schwartz function. We have
$$
(f * \delta_\Lambda)(\phi) =
(\widehat{f}\widehat{\delta_\Lambda}) (\ft\phi) = 
\widehat{\delta_\Lambda}(\ft\phi \ft f).
$$
However, this is
$$
\ft\phi(0) \widehat{f}(0) \widehat{\delta_\Lambda}(\Set{0}),
$$
as $\widehat{\delta_\Lambda}$, being a measure,
kills any continuous function vanishing on
$\supp \widehat{\delta_\Lambda}$.
Since $\phi$ is arbitrary we conclude that $f*\delta_\Lambda$ is a
constant.
\Qed

We shall need a different version of Theorem \ref{th:tiling-smooth}
here.
In the theorem that follows compact support and nonnegativity of $\ft f$
compensate for its lack of smoothness.
This theorem has also been proved and used by the author in \cite{nonsym}.
%%%%%%%%%%%%%%%%%%%%%%%%%%%%%%%%%%%%%%%%%%%%%%%%%%%%%%%%%%%%%%%%
%%%%%%%%%%%%%%% thm TILING IMPLIES SUPPORT
%%%%%%%%%%%%%%%%%%%%%%%%%%%%%%%%%%%%%%%%%%%%%%%%%%%%%%%%%%%%%%%%
\begin{theorem}\label{th:tiling-implies-supp}
Suppose that $f\ge 0$ is not identically $0$, that $f \in L^1(\RR^d)$,
$\widehat{f}\ge 0$ has compact support and $\Lambda\subset\RR^d$.
If $f+\Lambda$ is a tiling then
\beql{supp-cond-1}
\supp \ft{\delta_\Lambda} \subseteq \Set{\ft{f} = 0} \cup \Set{0}.
\eeq
\end{theorem}
{\bf Proof of Theorem \ref{th:tiling-implies-supp}.}
Assume that $f+\Lambda = w\RR^d$ and let
$$
K = \Set{\widehat{f} = 0} \cup \Set{0}.
$$
We have to show that
$$
\widehat{\delta_\Lambda}(\phi) = 0,
\ \ \forall \phi\in C_c^\infty(K^c).
$$
Since $\ft{\delta_\Lambda}(\phi) = \delta_\Lambda(\ft\phi)$
this is equivalent to $\sum_{\lambda\in\Lambda}\widehat\phi(\lambda) = 0$, for
each such $\phi$.
Notice that $h = \phi / \widehat{f}$ is a continuous function, but
not necessarily smooth. We shall need that $\widehat h \in L^1$.
This is a consequence of a well-known theorem of Wiener \cite[Ch.\ 11]{R73}.
We denote by $\TT^d = \RR^d/\ZZ^d$ the $d$-dimensional torus.
\begin{thm} (Wiener)\\
If $g \in C(\TT^d)$ has an absolutely convergent Fourier series
$$
g(x) = \sum_{n\in\ZZ^d} \widehat g(n) e^{2\pi i \Inner{n}{x}},
\ \ \ \widehat g \in \ell^1(\ZZ^d),
$$
and if $g$ does not vanish anywhere on $\TT^d$ then
$1/g$ also has an absolutely convergent Fourier series.
\end{thm}
Assume that
$$
\supp\phi,\ \supp\widehat{f} \subseteq \left(-{L\over2},{L\over2}\right)^d.
$$
Define the function $F$ to be:\\
(i) periodic in $\RR^d$ with period lattice $(L\ZZ)^d$,\\
(ii) to agree with $\widehat{f}$ on $\supp\phi$,\\
(iii)to be non-zero everywhere and,\\
(iv) to have $\widehat{F} \in \ell^1(\ZZ^d)$, i.e.,
$$
\widehat F = \sum_{n\in\ZZ^d} \widehat F(n) \delta_{L^{-1} n},
$$
is a finite measure in $\RR^d$.

One way to define such an $F$ is as follows.
First, define the $(L\ZZ)^d$-periodic function $g\ge 0$
to be $\ft{f}$ periodically extended.
The Fourier coefficients of $g$ are $\ft{g}(n) = L^{-d} f(-n/L) \ge 0$.
Since $g, \ft g \ge 0$ and $g$ is continuous at $0$ it is
easy to prove that $\sum_{n\in\ZZ^d} \ft g (n) = g(0)$, and therefore
that $g$ has an absolutely convergent Fourier series.
%By the Poisson Summation Formula (PSF) for $f$
%$$
%\sum_{n\in\ZZ^d} \ft{g}(n) = \sum_{n\in\ZZ^d} L^{-d} f\pp{n\over L}
% = \sum_{n\in\ZZ^d} \ft{f}(L n) < \infty
%$$
%since $\ft{f}$ has compact support.
%We should note here that there is no problem for the application
%of the PSF since both $f$ and $\ft f$ are nonnegative. This permits
%us to apply the PSF to the pair $\phi = \psi_\epsilon*f$,
%$\widehat\phi = \widehat{\psi_\epsilon}\widehat{f}$,
%where $\psi_\epsilon$ is a smooth, positive and positive-definite approximate
%identity. Then both $\phi$ and $\widehat\phi$ are Schwartz functions
%and the PSF is applicable.
%In the limit we get the PSF for $f$ and $\widehat{f}$.

Let $\epsilon$ be small enough to guarantee that $\ft{f}$
(and hence $g$) does not vanish on $(\supp\phi) + B_\epsilon(0)$.
Let $k$ be a smooth $(L\ZZ)^d$-periodic function which is equal to $1$
on $(\supp\phi)+(L\ZZ^d)$ and equal to $0$ off
$(\supp\phi + B_\epsilon(0))+(L\ZZ^d)$, and
satisfies $0\le k \le 1$ everywhere.
Finally, define
$$
F = k g + (1-k).
$$
Since both $k$ and $g$ have absolutely summable Fourier series and this
property is preserved under both sums and products, it follows that $F$ also
has an absolutely summable Fourier series. And by the nonnegativity of $g$ we
get that $F$ is never $0$, since $k=0$ on $Z(\ft f)+(L\ZZ^d)$.

By Wiener's theorem, $\widehat{F^{-1}} \in \ell^1(\ZZ^d)$, i.e., $\ft{F^{-
1}}$ is a finite
measure on $\RR^d$.
We now have that
$$
\FT{{\phi \over \widehat f}} = \widehat{\phi F^{-1}} =
  \widehat\phi * \widehat{F^{-1}} \in L^1(\RR^d).
$$
This justifies the interchange of the summation and integration below:
\begin{eqnarray*}
\sum_{\lambda\in\Lambda} \widehat\phi(\lambda) 
 &=& \sum_{\lambda\in\Lambda} \FT{{\phi\over\widehat{f}} \widehat{f}}
(\lambda) \\
 &=& \sum_{\lambda\in\Lambda} \FT{{\phi\over\widehat f}} *
\widehat{\widehat{f}}~(\lambda) \\
 &=& \sum_{\lambda\in\Lambda} \int_{\RR^d} \FT{{\phi\over\widehat f}}(y) f(y-
\lambda) ~dy \\
 &=& \int_{\RR^d} \FT{{\phi\over\widehat f}}(y)
  \sum_{\lambda\in\Lambda} f(y-\lambda) ~dy\\
 &=& w \int_{\RR^d} \FT{{\phi\over\widehat f}}(y)~dy\\
 &=& w{\phi\over\ft{f}}(0)\\
 &=& 0,
\end{eqnarray*}
as we had to show.
\Qed

In the other direction assume that we have
\beql{the-usual-condition}
\supp \ft{\delta_\Lambda} \subseteq \Set{\ft f = 0} \cup \Set{0}
\eeq
for some non-zero $f\ge 0$ in $L^1$ and that $\Lambda$ is of bounded density.
Since $\ft f(0) = \int f > 0$ it follows that in some neighborhood $N$ of $0$
we have $(\supp\ft{\delta_\Lambda}) \cap N = \Set{0}$.
Hence the set
\beql{O-set}
O = \left( \supp \ft{\delta_\Lambda} \setminus \Set{0} \right)^c
\eeq
is open and
$$
\Set{\ft f \neq 0} \subseteq O.
$$
%%%%%%%%%%%%%%%%%%%%%%%%%%%%%%%%%%%%%%%%%%%%%%%%%%%%%%%%%%%%%%%%
%%%%%%%%%%%%%%% thm SUPPORT IMPLIES TILING
%%%%%%%%%%%%%%%%%%%%%%%%%%%%%%%%%%%%%%%%%%%%%%%%%%%%%%%%%%%%%%%%
\begin{theorem}\label{th:supp-implies-tiling}
Suppose that $0\le f \in L^1(\RR^d)$, $\int f =1$,
$\Lambda$ (of uniformly bounded density) is of density $1$,
and that \eqref{the-usual-condition} holds.
Suppose also that for the open set $O$ of \eqref{O-set} and for each
$\epsilon>0$ there exists $f_\epsilon \ge 0$ in $L^1(\RR^d)$ such that
$\ft{f_\epsilon}$ is in $C^\infty$, $\supp \ft{f_\epsilon} \subseteq O$ and
$$
\Lone{f - f_\epsilon} \le \epsilon.
$$
Then $f + \Lambda$ is a tiling.
\end{theorem}

\noindent
{\bf Example.}
All bounded open convex sets $O$ have the property
required by the theorem, for all functions $f\ge 0$ such that
$\ft f$ is non-zero only in $O$.
To see this, assume, without loss of
generality, that $0 \in O$.
Construct then the functions ($\epsilon \to 0$)
$$
\ft{f_\epsilon} (x) = \psi_\epsilon(x) *
   \ft{f}\left({x\over 1-\epsilon}\right),
$$
where $\psi_\epsilon$ is a smooth, positive-definite approximate
identity supported in $(\epsilon/2) O$.
Then $\ft{f_\epsilon}$ is smooth, supported properly in
$O$ and $f_\epsilon$ converges to $f$ in $L^1$
(for example, by the dominated convergence theorem).
We will not need the above observation about convex domains below.

\noindent
{\bf Proof of Theorem \ref{th:supp-implies-tiling}.}
Suppose that $f_\epsilon$ is as in the Theorem.
First we show that $(\int f_\epsilon)^{-1}f_\epsilon + \Lambda$
is a tiling.
That is, we show that the convolution $f_\epsilon * \delta_\Lambda$
is a constant.
Let $\phi$ be $C_c^\infty$ function.
Then
$$
(f_\epsilon * \delta_\Lambda)(\phi) = 
\ft{f_\epsilon}\ft{\delta_\Lambda}(\ft\phi) =
\ft{\delta_\Lambda}(\ft\phi \ft{f_\epsilon}).
$$
But the function $\ft\psi = \ft\phi \ft{f_\epsilon}$ is a $C_c^\infty$
function whose support intersects $\supp\ft{\delta_\Lambda}$ only at $0$.
And, it is not hard to show, because $\Lambda$ has density $1$,
that $\ft{\delta_\Lambda}$ is equal to $\delta_0$ in a
neighborhood of $0$ (see \cite{polygons}).
Hence 
$$
(f_\epsilon * \delta_\Lambda)(\phi) =
 \left(\ft\phi \ft{f_\epsilon}\right) (0) =
 \int \phi \int f_\epsilon,
$$
and, since this is true for an arbitrary $C_c^\infty$ function $\phi$,
we conclude that $f_\epsilon * \delta_\Lambda = \int f_\epsilon$,
as we had to show.

For any set $\Lambda$ of uniformly bounded density we have ($B$ is any ball
in $\RR^d$ and $g\in L^1(\RR^d)$)
$$
\int_B\Abs{\sum_{\lambda\in\Lambda} g(x-\lambda)} ~dx \le C_{B,\Lambda}
	\int_{\RR^d}\Abs{g},
$$
(See \cite{KL96} for a proof of this in dimension $1$, which holds
for any dimension.)
Applying this for $g=f-f_\epsilon$ we obtain that
$$
\sum_{\lambda\in\Lambda} f_\epsilon(x-\lambda) \to
 \sum_{\lambda\in\Lambda} f(x-\lambda), \ \ \ \mbox{in $L^1(B)$}.
$$
Since $B$ is arbitrary this
implies that $\sum_{\lambda\in\Lambda} f(x-\lambda) = 1$, a.e.\ in $\RR^d$.
\Qed

We write $\widetilde f(x) = \overline{f(-x)}$.

Let $\Omega \subset \RR^d$ be a bounded open set of measure $1$,
$\one_\Omega$ its indicator function and
$f$ be such that $\ft{f} = \one_\Omega * \widetilde{\one_\Omega}$.
Then $\widetilde{f} = \Abs{\ft{\one_\Omega}}^2 \ge 0$,
$\int f = 1$ by Parseval's theorem.
Clearly we have $\Set{\ft{f} \neq 0} = \Omega - \Omega$.

Write
$$
\Omega_\epsilon = \Set{x\in\Omega:\ {\rm dist\,}(x, \partial\Omega) >
\epsilon},
$$
and define $f_\epsilon$ by
$$
\ft{f_\epsilon} = \psi_\epsilon * \one_{\Omega_\epsilon} *
 (\psi_\epsilon * \one_{\Omega_\epsilon})^{\widetilde{\ }}
$$
(or $\widetilde{f_\epsilon} = \Abs{\ft{\psi_\epsilon}}^2
 \Abs{\ft{\one_{\Omega_\epsilon}}}^2$),
where $\psi_\epsilon$ is a smooth, positive-definite
approximate identity supported in $B_{\epsilon/2}(0)$.

One can easily prove the following proposition.

\noindent
{\em If $g_n \to g$ in $L^2$ then $\Abs{g_n}^2 \to \Abs{g}^2$ in $L^1$.}

(For the proof just notice the identity
$$
\Abs{g}^2 - \Abs{g_n}^2 = \Abs{g - g_n}^2 +
  2\cdot{\rm Re} \left(\overline{g_n}(g-g_n)\right),
$$
integrate and use the triangle and Cauchy-Schwartz inequalities.)

Since $\psi_\epsilon*\one_{\Omega_\epsilon} \to \one_\Omega$ in $L^2$
(dominated convergence) we have (Parseval)
that $\ft{\psi_\epsilon}\ft{\one_{\Omega_\epsilon}} \to \ft{\one_\Omega}$
in $L^2$ and, using the proposition above, that
$\Abs{\ft{\psi_\epsilon}}^2 \Abs{\ft{\one_{\Omega_\epsilon}}}^2
  \to \Abs{\ft{\one_\Omega}}^2$ in $L^1$,
which means that $f_\epsilon \to f$ in $L^1$.

We also have that
$$
\supp{\ft{f_\epsilon}} \subseteq
 \overline{\Omega_{\epsilon/2}} - \overline{\Omega_{\epsilon/2}}
\subseteq \Omega - \Omega = \Set{\ft f \neq 0}.
$$

The assumptions of Theorem \ref{th:supp-implies-tiling} are therefore
satisfied.
Combining Theorems \ref{th:tiling-implies-supp} and
\ref{th:supp-implies-tiling} with the above observations we
obtain the following characterization of tiling which we will use
throughout the rest of the paper.
%%%%%%%%%%%%%%%%%%%%%%%%%%%%%%%%%%%%%%%%%%%%%%%%%%%%%%%%%%%%%%%%
%%%%%%%%%%%%%%% thm SUPPORT EQUIVALENT TO TILING
%%%%%%%%%%%%%%%%%%%%%%%%%%%%%%%%%%%%%%%%%%%%%%%%%%%%%%%%%%%%%%%%
\begin{theorem}\label{th:supp-equiv-tiling}
Let $\Omega$ be a bounded open set, $\Lambda$ a discrete set in $\RR^d$,
and $\delta_\Lambda = \sum_{\lambda\in\Lambda} \delta_\lambda$.
Then $\Abs{\ft{\one_\Omega}}^2 + \Lambda$ is a tiling if and only if
$\Lambda$ has uniformly bounded density and 
$$
(\Omega - \Omega) \cap \supp\ft{\delta_\Lambda} = \Set{0}.
$$
\end{theorem}

\mysection{Size of orthogonal packing regions. Spectra of the cube.}
\label{sec:size-ortho-regions}

The following theorem was conjectured in \cite{LRW98} (Conjecture 2.1).
%%%%%%%%%%%%%%%%%%%%%%%%%%%%%%%%%%%%%%%%%%%%%%%%%%%%%%%%%%%%%%%%
%%%%%%%%%%%%%%% thm SIZE OF ORTHOGONAL PACKING REGIONS
%%%%%%%%%%%%%%%%%%%%%%%%%%%%%%%%%%%%%%%%%%%%%%%%%%%%%%%%%%%%%%%%
\begin{theorem}\label{th:LRW-conj}
If $\Omega$ has measure $1$ and tiles $\RR^d$ then $\Abs{D} \le 1$
for any orthogonal packing region $D$ of $\Omega$.
\end{theorem}
{\bf Proof of Theorem \ref{th:LRW-conj}.}
Assume that $\Omega + \Lambda$ is a tiling.
Then $\dens \Lambda = 1$.
By Theorem \ref{th:tiling-smooth}
$$
\supp\ft{\delta_\Lambda} \subseteq Z(\ft{\one_\Omega}) \cup \Set{0}.
$$
Since $D$ is an orthogonal packing region for $\Omega$ we have
by definition,
and since $D-D = \Set{\one_D*\widetilde{\one_D} \neq 0}$,
$$
Z(\ft{\one_\Omega}) \subseteq Z(\one_D * \widetilde{\one_D}).
$$
Therefore
$$
\supp\ft{\delta_\Lambda} \subseteq Z(\one_D * \widetilde{\one_D})
	\cup \Set{0},
$$
and by Theorem \ref{th:supp-equiv-tiling} we obtain that
$\Abs{\ft{\one_D}}^2 + \Lambda$ is a tiling at level
$$
\dens\Lambda \int \Abs{\ft{\one_D}}^2 = \Abs{D}\ \ \mbox{(Parseval)}.
$$
On the other hand, the level of the tiling is (evaluating at $0$)
$$
\sum_\lambda \Abs{\ft{\one_D}}^2(-\lambda) \ge \Abs{\ft{\one_D}}^2 (0) =
	\int \one_D * \widetilde{\one_D} = \Abs{D}^2.
$$
Hence $\Abs{D} \ge \Abs{D}^2$ or $\Abs{D} \le 1$.
\Qed

\begin{dfn}\label{def:tight}
(Tight orthogonal packing regions, tight spectral pairs)\\
The open set $D$ is called a \underline{tight orthogonal packing region}
for $\Omega$ if it is an orthogonal
packing region for $\Omega$ and $\Abs{D}=1$.

A pair $(\Omega, D)$ of bounded open sets in $\RR^d$ is called
a \underline{tight spectral pair} if each is a tight orthogonal packing region
for the other.
\end{dfn}

The following result has also been proved in \cite{LRW98} (Theorem 3.1)
but for a smaller class of admissible open sets $\Omega$.
\begin{theorem}\label{th:th-3.1}
Suppose $\Omega, \Omega'$ are bounded open sets in $\RR^d$ of measure $1$.
Suppose also that $\Lambda$ is an orthogonal set of exponentials for $\Omega$
and that $\Omega' + \Lambda$ is a packing.

Then $\Lambda$ is a spectrum for $\Omega$ if and only if
$\Omega' + \Lambda$ is a tiling.
\end{theorem}
{\bf Proof.}
Since $E_\Lambda$ is an orthogonal set in $L^2(\Omega)$ it follows
that $\Abs{\ft{\one_\Omega}}^2+\Lambda$ is a packing and
$\Lambda$ is a spectrum for $\Omega$ if and only if
$\Abs{\ft{\one_\Omega}}^2 + \Lambda$ is also a tiling.
Notice that
$$
\int \Abs{\ft{\one_\Omega}}^2 = \int \one_{\Omega'} = 1.
$$
By Theorem \ref{th:tiling-from-packing} it follows that
$\Omega' + \Lambda$ is a tiling if and only if
$\Abs{\ft{\one_\Omega}}^2 + \Lambda$ is a tiling,
as we had to show.
\Qed

The following theorem relates, for a tight spectral
pair $(\Omega, D)$, the tilings of $D$ with the spectra
of $\Omega$.
\begin{theorem}\label{th:pair}
Assume that $(\Omega, D)$ is a tight spectral pair. Then
$\Lambda$ is a spectrum of $\Omega$ if and only if 
$D + \Lambda$ is a tiling.
\end{theorem}
{\bf Proof.}
$\Lambda$ is a spectrum of $\Omega$ if and only if
$\Abs{\ft{\one_\Omega}}^2 + \Lambda$ is a tiling (Theorem
\ref{th:connections}).

``$\Longleftarrow$'' If $D + \Lambda$ is a tiling then
$\supp\ft{\delta_\Lambda} \subseteq Z(\ft{\one_D}) \cup \Set{0}$
(by Theorem \ref{th:tiling-smooth}),
which is a subset of $Z(\one_\Omega * \widetilde{\one_\Omega})
\cup\Set{0}$ (this is because $\Omega$ is an
orthogonal packing region for $D$).
Hence $\Abs{\ft{\one_\Omega}}^2 + \Lambda$ is a tiling
(Theorem \ref{th:supp-equiv-tiling}).

``$\Longrightarrow$'' In the other direction, suppose that
$\Abs{\ft{\one_\Omega}}^2 + \Lambda$ is a tiling.
Then $\Lambda$ is an orthogonal set for $\Omega$ and hence
$D + \Lambda$ is a packing because $D$ is an orthogonal
packing region for $\Omega$ (Theorem \ref{th:connections}(3)).
But Theorem \ref{th:tiling-from-packing}
implies then that $D + \Lambda$ is a tiling as well.
\Qed

Let $Q = (-1/2, 1/2)^d$ be the $0$-centered unit cube in $\RR^d$.
An easy calculation gives for the Fourier transform of $\one_Q$:
\beql{qft}
\ft{\one_Q}(\xi_1,\ldots,\xi_d) =
 \prod_{j=1}^d{\sin \pi \xi_j \over \pi \xi_j}.
\eeq
It follows that $\ft{\one_Q}$ vanishes precisely at those points
with at least one non-zero integer coordinate.
From the definition of an orthogonal packing region
it follows that $(Q, Q)$ is a tight spectral pair.
Hence we have the following, which has already been proved in
\cite{LRW98,IP98}.
\begin{corollary}\label{cor:cube}
$\Lambda$ is a spectrum of $Q$ if and only if $Q + \Lambda$ is a tiling.
\end{corollary}
Let us also mention another case, in dimension $1$, when Theorem
\ref{th:pair} applies:
the set $\Omega = (0,1/2) \cup (1,3/2)$ has itself as a tight orthogonal
packing region, hence all its spectra are also tiling sets for $\Omega$
and vice versa (see \cite{LRW98}).

\mysection{Generalization of a theorem of Keller}
\label{sec:keller}

In \cite{LRW98,IP98} the following old result of Keller \cite{K} was used
in order to prove that the spectra of the unit cube $Q$
of $\RR^d$ are exactly those sets
$\Lambda$ such that $Q+\Lambda$ is tiling
(our Corollary \ref{cor:cube}).
\begin{thm} (Keller)\\
If $Q+\Lambda$ is a tiling and $0 \in \Lambda$ then
each non-zero $\lambda\in\Lambda$ has a non-zero integer coordinate.
\end{thm}
The use of Keller's theorem is avoided in our approach.
Furthermore it is an easy consequence of what we have already proved.

Having computed $\ft{\one_Q}$ and
its zero-set in the previous paragraph it is now evident
that Keller's theorem is a special case of the following result.
\begin{theorem}\label{th:generalized-keller}
Assume that $(\Omega,D)$ is a tight spectral pair and that
$\Omega+\Lambda$ is a tiling ($0 \in \Lambda$).
Then
$$
\Lambda \setminus \Set{0} \subseteq Z(\ft{\one_D}).
$$
\end{theorem}
{\bf Proof of Theorem \ref{th:generalized-keller}.}
Since $\Omega+\Lambda$ is a tiling we obtain that
$\Abs{\ft{\one_D}}^2 + \Lambda$ is a tiling (i.e., $\Lambda$ is a
spectral set for $D$) and, since $\Abs{\ft{\one_D}}^2(0) = 1$,
we obtain
$$
\ft{\one_D}(\lambda) = 0,\ \ (\lambda \in \Lambda\setminus\Set{0}).
$$
\Qed

\mysection{Bibliography}


\begin{thebibliography}{99}
\vspace{-1.5cm}

\bibitem{F74} B. Fuglede, Commuting self-adjoint partial differential
operators and a group theoretic problem,
J. Funct.\ Anal.\ {\bf 16}(1974), 101-121.

\bibitem{IKP99} A. Iosevich, N. Katz and S. Pedersen, Fourier bases and a
distance problem of Erd\H os, Math.\ Res.\ Letters {\bf 6} (1999), 251-255.

\bibitem{IP98} A. Iosevich and S. Pedersen, Spectral and
Tiling Properties of the Unit Cube, Inter.\ Math.\ Res.\ Notices {\bf 16},
819-828.

\bibitem{JP} P.E.T. Jorgensen and S. Pedersen, Spectral pairs in
cartesian coordinates, preprint.

%\bibitem{K68} Y. Katznelson, {\em An introduction to harmonic
%analysis}, Wiley (1968) and Dover (1976).

\bibitem{K} O.-H. Keller, \"Uber die l\"uckenlose Erf\"ullung des Raumes
mit W\"urfeln, J. Reine Angew.\ Math.\ {\bf 163} (1930), 231-248.

\bibitem{nonsym} M.N. Kolountzakis, Non-symmetric convex domains have
no basis of exponentials, Illinois J. Math., to appear.

\bibitem{polygons} M.N. Kolountzakis, On the structure of multiple
translational tilings by polygonal regions, preprint.

\bibitem{KL96} M.N. Kolountzakis and J.C. Lagarias,
Structure of tilings of the line by a function, Duke Math.\ J.

\bibitem{LRW98} J.C. Lagarias, J.A. Reeds and Y. Wang,
Orthonormal bases for exponentials for the $n$-cube, preprint (1998).

\bibitem{LS} J.C. Lagarias and P.W. Shor, Keller's
cube-tiling conjecture is false in high dimensions,
Bull.\ Amer.\ Math.\ Soc\ (N.S.) {\bf 27}(1992), 2, 279-283.

\bibitem{M80} P. McMullen, Convex bodies which tile space
by translation, Mathematika {\bf 27} (1980), 113-121.

\bibitem{R73} W. Rudin, {\em Functional analysis},
McGraw-Hill, New York, 1973.

\end{thebibliography}
\end{document}